\documentclass[11pt,reqno]{amsart}
\usepackage[T1]{fontenc}
\usepackage{lmodern}
\usepackage[a4paper,textwidth=6.5in,textheight=9in,hcentering]{geometry}
\usepackage{amsmath,amssymb,amsthm,mathtools,mathrsfs,bm}
\usepackage{booktabs,tabularx,array,enumitem,microtype}
\usepackage{xcolor}

\usepackage{hyperref}

\hypersetup{hidelinks,
 pdftitle={Entropy and frame measures for self-affine measures},
 pdfauthor={Zhiyong Wang and Wei Tang}}
\numberwithin{equation}{section}
\newtheorem{theorem}{Theorem}[section]
\newtheorem{lemma}[theorem]{Lemma}
\newtheorem{proposition}[theorem]{Proposition}
\newtheorem{corollary}[theorem]{Corollary}
\theoremstyle{definition}

\newtheorem{example}[theorem]{Example}
\theoremstyle{remark}

\newcommand{\R}{\mathbb R}
\newcommand{\Z}{\mathbb Z}
\newcommand{\N}{\mathbb N}

\newcommand{\E}{\mathbb E}
\newcommand{\PP}{\mathbb P}
\newcommand{\ee}{\mathrm e}
\newcommand{\ind}{\mathbf 1}
\DeclareMathOperator{\supp}{supp}

\setlist[enumerate,1]{label=\textup{(\roman*)},leftmargin=2.1em,itemsep=3pt,topsep=5pt}
\allowdisplaybreaks[1]

\title[Entropy and frame measures]{Entropy and frame measures for self-affine measures}
\author[Z. Wang]{Zhiyong Wang}
\address{School of Mathematics and Statistics, Hunan First Normal University, Changsha, Hunan 410205, P.~R.~China}
\email{wzyzzql@163.com}
\author[W. Tang]{Wei Tang$^*$}
\address{School of Mathematics and Statistics, Hunan First Normal University, Changsha, Hunan 410205, P.~R.~China}
\email{twmath2016@163.com}
\thanks{The first author is supported by the Natural Science Foundation of Hunan Province (Grant No.~2024JJ4017) and the Science and Technology Innovation Program of Hunan Province (Grant No.~2024RC3227). The second author is supported by the National Natural Science Foundation of China (Grant No.~12271156) and the Natural Science Foundation of Hunan Province (Grant No.~2023JJ40201). Corresponding author: Wei Tang.}
\date{November 11, 2026}

\begin{document}
\begin{abstract}
Let $M$ be a real expanding $d\times d$ matrix and let
$\mathscr D=\{d_0,\ldots,d_{N-1}\}\subset\R^d$ consist of
$N\ge2$ distinct vectors. We study the existence of frame measures
for the self-affine probability measure $\mu$ generated by the maps
$S_i(x)=M^{-1}(x+d_i)$ with equal weights $1/N$.
Extending the entropy obstruction for self-similar measures
with a common scalar contraction, we prove that $\mu$
admits no frame measure if the conditional entropy of
the first symbol given the spatial point is positive.
In particular, this holds whenever $N>|\det M|$, without assuming absolute continuity
or imposing arithmetic conditions.
When $M$ has integer entries, $\mathscr D\subset\Z^d$, and
$N=|\det M|$, the measure $\mu$ admits a frame measure if and only
if $\mu$ is normalized Lebesgue measure on its attractor.
\end{abstract}
\maketitle

\section{Introduction}

Let $\nu$ be a compactly supported Borel probability measure on $\R^d$.
For $f\in L^2(\nu)$, write
\[
 \widehat{f\nu}(\lambda)
 =\int_{\R^d}f(x)\ee^{-2\pi i\langle\lambda,x\rangle}\,d\nu(x),
 \qquad e_\lambda(x)=\ee^{2\pi i\langle\lambda,x\rangle}.
\]
The integral is well defined since $L^2(\nu)\subset L^1(\nu)$.
A positive Borel measure $\Gamma$ on $\R^d$ is called a
\emph{frame measure} for $\nu$ if there are constants
$0<A\le C<\infty$ such that
\begin{equation}\label{eq:frame}
 A\|f\|_{L^2(\nu)}^2
 \le\int_{\R^d}|\widehat{f\nu}(\lambda)|^2\,d\Gamma(\lambda)
 \le C\|f\|_{L^2(\nu)}^2
 \qquad(f\in L^2(\nu)).
\end{equation}
The upper inequality alone defines a \emph{Bessel measure}.
When $\Gamma=\sum_{\lambda\in\Lambda}\delta_\lambda$ for a countable set
$\Lambda\subset\R^d$, the inequalities in \eqref{eq:frame}
give a Fourier frame of exponentials. Positive
weights on the atoms give weighted Fourier frames. If the
exponentials form an orthonormal basis, then $\nu$ is called
\emph{spectral}, and $\Lambda$ is a \emph{spectrum}. Counting measure
on a spectrum is a frame measure with bounds $A=C=1$.
Frame measures were introduced by Dutkay, Han and Weber \cite{DHW}.

The existence of a frame measure imposes uniformity conditions on the underlying measure. Lai \cite[Theorem~1.1]{Lai11}
proved that a compactly supported absolutely continuous
probability measure admits a Fourier frame if and only if
its density is essentially bounded above and bounded away
from zero on the set where it is positive. Dutkay and Lai \cite[Theorems~2.2 and 2.4]{DL}
obtained uniformity conditions for measures admitting
frame measures. They also proved that, for an affine
IFS with a common expanding matrix and no overlap,
the invariant measure admits a frame measure only
if all probability weights are equal
\cite[Theorem~3.4]{DL}.
Without assuming no overlap, Deng and Chen
\cite[Theorem~1.1]{DC} proved that spectral self-affine measures
have equal probability weights and pairwise mutually singular
first-level component measures. These results motivate the study of frame measures for
overlapping self-affine measures.

In \cite[Theorem~1.1]{WT}, we proved that positive conditional
entropy of the first symbol given the spatial point rules
out frame measures for equal-weight self-similar measures
on $\R$. We also showed that this conditional entropy
vanishes if and only if the first-level component measures
are mutually singular. In particular, the entropy obstruction
excludes the entire supercritical range $1/N<\rho<1$ for
arbitrary sets of $N$ distinct real digits. The same paper
establishes a separate arithmetic obstruction for consecutive
digits in the range $0<\rho<1/N$.

W.-J. Wang \cite[Theorem~1.1]{Wang26} extended this obstruction
to equal-weight self-similar measures generated by
$S_i(x)=\rho x+b_i$ in $\R^d$. He proved that the existence
of a frame measure implies essential injectivity of the
coding map, vanishing conditional entropy at every finite
level, and exact dimensionality with Hausdorff dimension
$\log N/(-\log\rho)$. In this setting, without any frame assumption, these properties are equivalent to
$H(\omega_1\mid X)=0$ (see \eqref{eq:def-H}
for the definition of $H$).
We record the complete relation in
Proposition~\ref{prop:coding-rigidity} and
Corollary~\ref{cor:dimension-equivalence}. As noted in
\cite[Remark~4.8]{Wang26},
the one-dimensional Dirac conclusion follows from \cite{WT}.

Here we extend the entropy obstruction from scalar contractions $\rho$
to arbitrary real expanding matrices $M$. The main new ingredient is a geometric construction:
at each level, we find a measurable set on which only
one component measure contributes and whose mass is
a fixed positive proportion of that component's total mass
(see Lemmas~\ref{lem:uniform-gap}
and~\ref{lem:cap-restriction}). The selected branches need not be nested, unlike those
in the endpoint construction of \cite{WT} and the repeated
exposed-branch construction of \cite{Wang26}.
The remaining conditional-probability and Fourier-energy
arguments follow \cite{WT}.

Fix $d\ge1$ and $N\ge2$, and put $D_N=\{0,1,\ldots,N-1\}$.
Let $M$ be a real $d\times d$ expanding matrix; that is, every
eigenvalue of $M$ has modulus greater than one. Let
$\mathscr D=\{d_0,\ldots,d_{N-1}\}\subset\R^d$ consist of distinct
vectors. We denote by $\mu_{M,\mathscr D}$ the uniform self-affine
probability measure determined by
\begin{equation}\label{eq:general-ifs}
 \mu_{M,\mathscr D}
 =\frac1N\sum_{i=0}^{N-1}(S_i)_\#\mu_{M,\mathscr D},
 \qquad S_i(x)=M^{-1}(x+d_i).
\end{equation}
Here $(S_i)_\#\mu(B)=\mu(S_i^{-1}(B))$ denotes the pushforward.
An equivalent norm makes $M^{-1}$ a contraction, so existence and
uniqueness follow from \cite{Hutchinson}. On the uniform product
space $D_N^{\N}$, the corresponding random series is
\begin{equation}\label{eq:coding-intro}
 X=\sum_{j\ge1}M^{-j}d_{\omega_j}.
\end{equation}
The symbol $i\in D_N$ labels the actual digit $d_i\in\mathscr D$. We write $H(\omega_1\mid X)$ for the conditional entropy
of the first symbol $\omega_1$ given the spatial point $X$,
as defined in \eqref{eq:def-H}. It vanishes if and only if
$\omega_1$ is determined by $X$ almost surely.
Our first result gives an entropy obstruction to the
existence of a frame measure.

\begin{theorem}\label{thm:general-digits}
Let $\mu=\mu_{M,\mathscr D}$ be as in \eqref{eq:general-ifs}.
If $H(\omega_1\mid X)>0$, then $\mu$ admits no frame measure.
\end{theorem}

Consequently, the equivalent zero-entropy conditions in
Proposition~\ref{prop:coding-rigidity} are necessary for the
existence of a frame measure. In particular, the component measures are pairwise mutually
singular at every level, and the coding map is injective on
a full-measure subset of the symbolic space.
These properties do not require the branch images to be disjoint. The proof also constructs a sequence of unit-norm functions
in $L^2(\mu)$, depending only on $M$ and $\mathscr D$,
whose Fourier energies tend to zero for every Bessel measure
when $H(\omega_1\mid X)>0$.
The zero-entropy conditions are not sufficient for the existence
of a frame measure; see Example~\ref{ex:zero-not-frame}.

The projection entropy formula of Feng and Hu
\cite[Theorem~2.3]{FH}, together with an elementary counting
estimate, gives
\begin{equation*}
 H(\omega_1\mid X)\ge\log N-\log|\det M|.
\end{equation*}
We derive this estimate in Proposition~\ref{prop:determinant}.
Together with Theorem~\ref{thm:general-digits}, it gives
the following corollary.

\begin{corollary}\label{thm:overlap}
If $N>|\det M|$, then the measure $\mu_{M,\mathscr D}$ in
\eqref{eq:general-ifs} admits
no frame measure.
\end{corollary}

No lattice, separation, or absolute-continuity assumption is needed.
For $M=\rho^{-1}I_d$, the nonexistence condition reduces to
$N\rho^d>1$, as in \cite{Wang26}. More generally, the result applies to common similarities
with rotations (Corollary~\ref{cor:similarities})
and to nonconformal expanding matrices, including anisotropic
expansions and matrices with nontrivial Jordan blocks
(Example~\ref{ex:jordan}). The latter example has no exact overlaps. If all eigenvalues of $M$ have modulus $\rho^{-1}$,
Feng and Hu's dimension formula \cite[Theorem~2.8]{FH}
gives
$H(\omega_1\mid X)=0$ if and only if
$\dim_{\rm H}\mu=\log N/(-\log\rho)$,
even for matrices with nontrivial Jordan blocks.
For general expanding matrices, this dimension formula need
not hold with $\rho=|\det M|^{-1/d}$; see
Example~\ref{ex:anisotropic-dimension}.

We next consider the critical case $N=|\det M|$.
An \emph{exact overlap} occurs if $S_u=S_v$ for two distinct
words $u,v$ of the same length, where
$S_w=S_{w_1}\circ\cdots\circ S_{w_n}$.
For a Borel set $E\subset\R^d$, we write $|E|$
for its Lebesgue measure.

\begin{theorem}\label{thm:critical}
Suppose that $M$ is an expanding integer matrix,
$\mathscr D\subset\Z^d$, and $\#\mathscr D=N=|\det M|$.
Put $\mu=\mu_{M,\mathscr D}$ and $K=\supp\mu$.
The following conditions are equivalent.
\begin{enumerate}
\item The measure $\mu$ admits a frame measure.
\item The iterated function system has no exact overlap.
\item For every $n\ge1$,
$
 \#\bigl(\mathscr D+M\mathscr D+\cdots+M^{n-1}\mathscr D\bigr)=N^n.
$
\item $|K|>0$.
\item $|K|>0$ and $\mu=|K|^{-1}\ind_K\,dx$.
\item $H(\omega_1\mid X)=0$.
\end{enumerate}
Under these conditions, $d\Gamma(\lambda)=|K|\,d\lambda$ is a
Parseval frame measure for $\mu$.
\end{theorem}

The equivalence between the cardinality condition and positive
Lebesgue measure is classical in the theory of integral self-affine
tiles; see Lagarias and Wang
\cite[Theorem~1.1 and Corollary~1.1a]{LW}.
Theorem~\ref{thm:general-digits} and Lemma~\ref{lem:exact-overlap}
show that the existence of a frame measure rules out exact overlaps. Conversely, if there are no exact overlaps, the classical
result gives $|K|>0$, and Plancherel's theorem yields
a frame measure.
We give a proof, based on lattice counting, that the
cardinality condition implies
$|K|>0$. The digits need not form a complete set of
representatives for $\Z^d/M\Z^d$.

Section~\ref{sec:prelim} recalls the conditional-probability
identities from \cite{WT}, proves the coding equivalences,
and uses the results of Feng and Hu \cite{FH} to obtain
the dimension statements and determinant bound.
Section~\ref{sec:overlap} proves Theorem~\ref{thm:general-digits}
and Corollary~\ref{thm:overlap}. Section~\ref{sec:critical} proves
Theorem~\ref{thm:critical} and gives examples illustrating the
matrix and dimensional hypotheses.

\section{Preliminaries}\label{sec:prelim}

Throughout the paper, $\N=\{1,2,\ldots\}$, all logarithms are
natural, and we use the convention $0\log0=0$. We write $\mu=\mu_{M,\mathscr D}$ and
$K=\supp\mu$. For a Borel set $E$, we denote the restriction of a measure
$\theta$ to $E$ by $\theta|_E$.
The Hausdorff dimension of $\mu$ is defined by
\[
 \dim_{\rm H}\mu
 =\inf\{\dim_{\rm H}E:E\text{ is Borel and }\mu(E)=1\}.
\]
We say that $\mu$ is \emph{exact dimensional} if there exists
a constant $s\ge0$ such that
\[
\lim_{r\downarrow0}
\frac{\log\mu(B(x,r))}{\log r}=s
\quad\text{for }\mu\text{-almost every }x.
\]

\subsection{Coding and conditional probabilities}

Let $\Omega=D_N^{\N}$ with uniform product probability $\PP$, and
let $\vartheta$ be the left shift. Since $M$ is expanding, there exist constants $C>0$ and
$0<\tau<1$ such that $\|M^{-j}\|\le C\tau^j$, where $\|\cdot\|$ denotes the Euclidean operator norm.
Thus the series in \eqref{eq:coding-intro} converges uniformly.
The random variable $X$ defined by \eqref{eq:coding-intro}
has distribution $\mu$. Write $\pi(\omega)=X(\omega)$ for the continuous coding map, and set
$X_k=X\circ\vartheta^k$. For a word
$u=(u_1,\ldots,u_n)\in D_N^n$, put
\begin{equation}\label{eq:word-maps}
 t_u=\sum_{j=1}^n M^{-j}d_{u_j},\qquad
 b_u=M^nt_u,\qquad
 S_u(x)=M^{-n}x+t_u=M^{-n}(x+b_u).
\end{equation}
Thus $X=S_{\omega|_n}(X_n)$, where
$\omega|_n=(\omega_1,\ldots,\omega_n)$, and $X_n$ is independent
of $\omega|_n$ with distribution $\mu$.

For each finite word $u$, define the branch measure $\mu_u$
and its Radon--Nikodym derivative $p_u$ by
\begin{equation}\label{eq:def-pw}
\mu_u=N^{-|u|}(S_u)_\#\mu,\qquad
p_u=\frac{d\mu_u}{d\mu}.
\end{equation}
Since $\sum_{|u|=n}\mu_u=\mu$, we have $0\le\mu_u\le\mu$.
The derivatives therefore exist, and we choose Borel versions
taking values in $[0,1]$ and vanishing outside $K$.
As shown below, for $|u|=n$, $p_u(x)$ is a version of
$\mathbb P(\omega|_n=u\mid X=x)$ for $\mu$-almost every $x$,
without choosing a coding of $x$. Set
\begin{equation}\label{eq:def-H}
H:=H(\omega_1\mid X)
=-\int\sum_{i=0}^{N-1}p_i(x)\log p_i(x)\,d\mu(x),
\qquad
R_n=\int\sum_{|u|=n}p_u(x)^2\,d\mu(x),
\end{equation}
and
\[
H_n:=H(\omega_1,\ldots,\omega_n\mid X)
=-\int\sum_{|u|=n}p_u(x)\log p_u(x)\,d\mu(x).
\]
Thus $H_1=H$. We retain the notation $R_n$ from \cite{WT}.
The finite-level conditional
entropy $H_n$ is denoted by $R_n$ in \cite{Wang26}.

For a word $u$ of length $n$, let
$[u]=\{\omega\in\Omega:\omega|_n=u\}$ be the symbolic cylinder.
Since $\Omega$ and $K$ are compact metric spaces, there exists
a regular conditional distribution $\PP_x$ of $\omega$
given $X=x$. Thus, for Borel sets $A\subset\Omega$ and $B\subset K$,
\begin{equation}\label{eq:disintegration-formula}
\PP(A\cap\pi^{-1}(B))=\int_B\PP_x(A)\,d\mu(x),
\qquad
\PP_x(\pi^{-1}(\{x\}))=1\quad\mu\text{-a.e.}
\end{equation}

The following identities and decay result are the
multidimensional forms of \cite[Lemmas~3.1, 3.4 and 3.5]{WT}.
Their proofs use the Bernoulli coding and the invertibility
of the maps $S_i$.

\begin{lemma}\label{lem:prefix-probabilities}
For every $n\ge1$ and $u\in D_N^n$,
\[
p_u(X)=\PP(\omega|_n=u\mid X)\quad\PP\text{-a.s.},
\qquad 0\le p_u\le1,\qquad
\sum_{|v|=n}p_v=1\quad\mu\text{-a.e.}
\]
In particular, $\PP_x([u])=p_u(x)$ for all finite words
outside a common $\mu$-null set. Moreover,
\begin{equation}\label{eq:prefix-product}
P_n:=p_{\omega|_n}(X)
=\prod_{j=0}^{n-1}p_{\omega_{j+1}}(X_j)
\quad\PP\text{-a.s.},\qquad R_n=\E P_n.
\end{equation}
If $H>0$, then $R_n\to0$.	
\end{lemma}

\begin{proof}
The tail $X_n$ has distribution $\mu$ and is independent
of the first $n$ symbols. Since $X=S_u(X_n)$ on $[u]$,
for every Borel set $B$ we have
\[
\PP(X\in B,\ \omega|_n=u)
=N^{-n}\mu(S_u^{-1}(B))
=\mu_u(B)
=\int_B p_u(x)\,d\mu(x),
\]
where the last two equalities follow from \eqref{eq:def-pw}.
This proves the conditional-probability identity.
The disintegration formula \eqref{eq:disintegration-formula} also gives
$\PP_x([u])=p_u(x)$ for $\mu$-almost every $x$.
The bounds $0\le p_u\le1$ and the identity
$\sum_{|u|=n}p_u=1$ follow from
$0\le\mu_u\le\mu$ and $\sum_{|u|=n}\mu_u=\mu$.

For $i\in D_N$ and a nonempty word $v$, we have
$\mu_{iv}=N^{-1}(S_i)_\#\mu_v$. Thus, for every bounded
Borel function $\varphi$,
\[
\int\varphi\,d\mu_{iv}
=\frac1N\int\varphi(S_i y)p_v(y)\,d\mu(y)
=\int\varphi(x)p_v(S_i^{-1}x)p_i(x)\,d\mu(x).
\]
By uniqueness of Radon--Nikodym derivatives,
$p_{iv}(x)=p_i(x)p_v(S_i^{-1}x)$ for $\mu$-almost every $x$.

Since there are only countably many finite words, all the
preceding identities hold outside a common $\mu$-null set
$\mathcal N$.
Each $X_j$ has distribution $\mu$, so almost surely
$X_j\notin\mathcal N$ for every $j\ge0$.
Using $S_{\omega_{j+1}}^{-1}(X_j)=X_{j+1}$, we can therefore
iterate the chain rule to obtain the product formula
in \eqref{eq:prefix-product}.
Moreover, since $P_n=p_u(X)$ on $[u]$, conditioning on $X$
gives
\[
\E(P_n\mid X)
=\sum_{|u|=n}p_u(X)\PP(\omega|_n=u\mid X)
=\sum_{|u|=n}p_u(X)^2.
\]
Taking expectations yields $\E P_n=R_n$.

It remains to prove that $R_n\to0$ when $H>0$.
Set $\Psi=1-p_{\omega_1}(X)$. Then $0\le\Psi\le1$ and
\[
\kappa:=\E\Psi
=1-R_1
=\int\left(1-\sum_{i=0}^{N-1}p_i(x)^2\right)\,d\mu(x).
\]
A probability vector has positive entropy exactly when
at least two entries are positive. Hence $H>0$ implies
that $\sum_i p_i(x)^2<1$ on a set of positive
$\mu$-measure, and therefore $\kappa>0$.
Since the left shift $\vartheta$ is ergodic, Birkhoff's ergodic theorem gives
\[
\frac1n\sum_{j=0}^{n-1}\Psi\circ\vartheta^j
\longrightarrow\kappa>0\quad\PP\text{-a.s.}
\]
Thus the sum tends to infinity almost surely.
Since $\Psi\circ\vartheta^j=1-p_{\omega_{j+1}}(X_j)$,
the product formula and $1-u\le\ee^{-u}$ give
\[
0\le P_n
\le\exp\left(-\sum_{j=0}^{n-1}\Psi\circ\vartheta^j\right)
\longrightarrow0\quad\PP\text{-a.s.}
\]
Finally, $0\le P_n\le1$, so dominated convergence gives
$R_n=\E P_n\to0$.	
\end{proof}

\subsection{The equivalent zero-entropy conditions}

The next proposition gives the identity $H_n=nH$ and
characterizes zero conditional entropy in terms of the
conditional distributions, the coding map, and the branch
measures. These statements do not require $\mu$ to admit a frame measure. Theorem~\ref{thm:general-digits} shows that the equivalent conditions are necessary for the existence of a frame measure. The equivalence between $H=0$ and mutual
singularity of the first-level branch measures is the criterion
in \cite[Lemma~3.6]{WT}. The remaining equivalences follow from
Lemma~\ref{lem:prefix-probabilities} and disintegration.
For the corresponding coding results in the homothetic case,
see \cite[Proposition~4.3 and Corollaries~4.4 and 4.5]{Wang26}.

\begin{proposition}\label{prop:coding-rigidity}
For every $n\ge1$, $H_n=nH$. Without any frame-measure assumption,
the following conditions are equivalent:
\begin{enumerate}
\item $H=0$.
\item $H_n=0$ for every $n\ge1$.
\item $\PP_x$ is a Dirac measure for $\mu$-almost every $x$.
\item There are Borel sets $\Omega_0\subset\Omega$ and $K_0\subset K$
with $\PP(\Omega_0)=\mu(K_0)=1$ such that
$\pi:\Omega_0\to K_0$ is a Borel bijection.
\item For every $n\ge1$, the measures $\mu_w$, $|w|=n$, are
pairwise mutually singular.
\end{enumerate}
In \textup{(ii)} and \textup{(v)}, ``for every $n\ge1$''
may be replaced by ``for some $n\ge1$''.
\end{proposition}

\begin{proof}
Since $P_n=p_u(X)$ on $[u]$, Lemma~\ref{lem:prefix-probabilities}
gives, for every $n\ge1$,
\[
\begin{aligned}
\PP(P_n=0)
=\sum_{|u|=n}
\PP(\omega|_n=u,\ p_u(X)=0)
=\sum_{|u|=n}
\int_{\{p_u=0\}}p_u\,d\mu=0.
\end{aligned}
\]
Conditioning on $X$ also gives
$\E[-\log P_n]=H_n\le n\log N$.
We may therefore take logarithms and expectations in
\eqref{eq:prefix-product}. Since
$p_{\omega_{j+1}}(X_j)=p_{\omega_1}(X)\circ\vartheta^j$
and $\PP$ is invariant under $\vartheta$, we obtain
\begin{equation}\label{eq:entropy-chain}
H_n=\E[-\log P_n]
=\sum_{j=0}^{n-1}\E[-\log p_{\omega_{j+1}}(X_j)]
=nH.
\end{equation}
This proves the equivalence of \textup{(i)} and \textup{(ii)}.

Next, fix $n\ge1$. For distinct words $w,v\in D_N^n$,
the relation $\mu_w\perp\mu_v$ is equivalent to
$p_wp_v=0$ $\mu$-almost everywhere.
Indeed, if $\mu_w\perp\mu_v$, there is a Borel set $A$ such that
$\mu_w(A^c)=\mu_v(A)=0$. Thus $p_w=0$ almost everywhere on
$A^c$ and $p_v=0$ almost everywhere on $A$, so $p_wp_v=0$
almost everywhere. Conversely, if $p_wp_v=0$ almost everywhere,
then $\{p_w>0\}$ has full $\mu_w$-measure and zero
$\mu_v$-measure. Since $\sum_{|w|=n}p_w=1$ almost everywhere, pairwise mutual
singularity at level $n$ is therefore equivalent to the following
condition: for almost every $x$, exactly one of the numbers
$p_w(x)$ equals $1$, and all the others equal $0$.
This is also equivalent to $H_n=0$, since $-t\log t$ is
nonnegative on $[0,1]$ and vanishes only at $0$ and $1$.
In this case, $\mu_w=\mu|_{\{p_w=1\}}$ for every $|w|=n$.
This proves the equivalence of \textup{(ii)} and \textup{(v)}.
By \eqref{eq:entropy-chain}, either condition may be checked
at a single level $n\ge1$.

We now prove that \textup{(ii)} implies \textup{(iii)}.
Since $H_n=0$ for every $n\ge1$ and there are only
countably many finite words, we can choose a Borel set $K_0\subset K$ with
$\mu(K_0)=1$ such that
$
\PP_x([u])=p_u(x)
$
for every $x\in K_0$ and every finite word $u$, and each probability vector
$(p_u(x))_{|u|=n}$ is a point mass. By shrinking $K_0$ if necessary, we may also assume that
$\PP_x(\pi^{-1}(\{x\}))=1$ for every $x\in K_0$.
Fix $x\in K_0$. For each $n$, let $u^{(n)}$ be the unique word
of length $n$ such that $\PP_x([u^{(n)}])=1$.
The length-$n$ prefix of $u^{(n+1)}$ also has conditional
probability one, so it must equal $u^{(n)}$.
Hence the cylinders $[u^{(n)}]$ are nested and their
intersection consists of a single sequence, denoted by $\omega(x)$.
Continuity from above gives
$\PP_x(\{\omega(x)\})=1$, so $\PP_x=\delta_{\omega(x)}$.
Moreover, the support property gives $\pi(\omega(x))=x$.
This proves \textup{(iii)}.
Conversely, if \textup{(iii)} holds, then
$\PP_x([u])=p_u(x)$ shows that every finite-level probability
vector is a point mass almost everywhere. Thus $H_n=0$ for
every $n$, proving \textup{(ii)}.

To prove \textup{(iv)} from \textup{(ii)}, let $K_0$ and
$\omega(x)$ be as constructed above, and set
\[
\Omega_0
=\{\omega\in\Omega:\pi(\omega)\in K_0,\
p_{\omega|_n}(\pi(\omega))=1
\text{ for every }n\ge1\}.
\]
For each $n\ge1$, let
\[
A_n
=\{\omega\in\Omega:p_{\omega|_n}(\pi(\omega))=1\}
=\bigcup_{|u|=n}
\bigl([u]\cap\pi^{-1}(\{p_u=1\})\bigr).
\]
Since $\pi$ is continuous and each $p_u$ is Borel,
the sets $A_n$ are Borel. Hence
$\Omega_0=\pi^{-1}(K_0)\cap\bigcap_{n\ge1}A_n$
is Borel. For every $x\in K_0$, we have $\omega(x)\in\Omega_0$
and $\pi(\omega(x))=x$.
If another sequence in $\Omega_0$ has image $x$,
its prefix of length $n$ must equal $u^{(n)}$ for
every $n\ge1$, and hence it equals $\omega(x)$.
Thus $\pi|_{\Omega_0}:\Omega_0\to K_0$ is a
Borel bijection, since $\pi$ is continuous.
Its inverse is also Borel, because the cylinders
generate the Borel $\sigma$-algebra of $\Omega$ and
$
\{x\in K_0:\omega(x)\in[u]\}
=\{x\in K_0:p_u(x)=1\}
$
is Borel for every finite word $u$.
Finally, $\PP_x=\delta_{\omega(x)}$ and
$\omega(x)\in\Omega_0$ for every $x\in K_0$.
By \eqref{eq:disintegration-formula},
\[
\PP(\Omega_0)
=\int_{K_0}\PP_x(\Omega_0)\,d\mu(x)
=\mu(K_0)=1.
\]

Conversely, assume \textup{(iv)}. By \eqref{eq:disintegration-formula}, we have
$
1=\PP(\Omega_0)
=\int_K\PP_x(\Omega_0)\,d\mu(x),
$
so $\PP_x(\Omega_0)=1$ for $\mu$-almost every $x$.
Together with
$\PP_x(\pi^{-1}(\{x\}))=1$ almost everywhere,
this gives
$
\PP_x(\Omega_0\cap\pi^{-1}(\{x\}))=1
$
for $\mu$-almost every $x\in K_0$.
Since $\pi|_{\Omega_0}$ is a bijection onto $K_0$,
the set $\Omega_0\cap\pi^{-1}(\{x\})$
consists of the single sequence
$\omega(x)=(\pi|_{\Omega_0})^{-1}(x)$.
Hence $\PP_x=\delta_{\omega(x)}$ for
$\mu$-almost every $x$, proving \textup{(iii)}.
\end{proof}

Condition \textup{(iv)} concerns a restriction to full-measure
sets; points in $K_0$ may have other codings outside $\Omega_0$.
Condition \textup{(v)} concerns branch measures at a fixed level;
the closed sets $S_w(K)$ may intersect.

\subsection{Projection entropy and dimension}

We recall projection entropy to relate the zero-entropy
conditions in Proposition~\ref{prop:coding-rigidity}
to the dimension of $\mu$.
For a finite Borel partition $\mathcal A$ of $\Omega$ and
a sub-$\sigma$-algebra $\mathcal G$ of $\mathcal B(\Omega)$,
write
\[
H_\PP(\mathcal A\mid\mathcal G)
=-\int\sum_{A\in\mathcal A}
\E[\ind_A\mid\mathcal G]
\log\E[\ind_A\mid\mathcal G]\,d\PP.
\]
Let $\mathcal P=\{[i]:i\in D_N\}$ and
$\mathcal G=\pi^{-1}\mathcal B(\R^d)$, where $[i]=\{\omega\in\Omega:\omega_1=i\}$.
The projection entropy of Feng and Hu
\cite[Definition~2.1]{FH} is
\begin{equation}\label{eq:projection-entropy}
h_\pi(\vartheta,\PP)
=H_\PP(\mathcal P\mid\vartheta^{-1}\mathcal G)
-H_\PP(\mathcal P\mid\mathcal G)
=\log N-H.
\end{equation}
Indeed, $\mathcal G=\sigma(X)$ and
$\vartheta^{-1}\mathcal G=\sigma(X_1)$.
Lemma~\ref{lem:prefix-probabilities} gives
$H_\PP(\mathcal P\mid\mathcal G)=H$.
Since $\omega_1$ is uniformly distributed on $D_N$
and is independent of $X_1$, we have
$H_\PP(\mathcal P\mid\vartheta^{-1}\mathcal G)=\log N$.
For the $k$-step system indexed by all $N^k$ words of
length $k$, the first block $(\omega_1,\ldots,\omega_k)$
is uniformly distributed and independent of $X_k$.
The same calculation and \eqref{eq:entropy-chain}
therefore give projection entropy
$h_\pi(\vartheta^k,\PP)
=k\log N-H_k
=k h_\pi(\vartheta,\PP)$.

The following dimension formula follows from
\cite[Theorem~2.8]{FH}, applied to a contractive iterate
if necessary. The hypothesis that all eigenvalues of
$M$ have the same modulus includes the case
$M^{-1}=\rho O$, where $0<\rho<1$ and $O$ is orthogonal,
but does not require $M$ to be diagonalizable.

\begin{corollary}\label{cor:dimension-equivalence}
Suppose that all eigenvalues of $M$ have the same modulus
$\rho^{-1}>1$. Then $\mu$ is exact dimensional and
\begin{equation}\label{eq:similarity-dimension}
\dim_{\rm H}\mu=\frac{\log N-H}{-\log\rho}.
\end{equation}
Consequently, each condition in
Proposition~\ref{prop:coding-rigidity} is equivalent to
$
\dim_{\rm H}\mu=\frac{\log N}{-\log\rho}.
$
These conclusions do not require a frame measure.
\end{corollary}

\begin{proof}
For an invertible matrix $A$, write
$s_{\min}(A)=\min_{\|v\|=1}\|Av\|=\|A^{-1}\|^{-1}$.
By the spectral-radius formula,
\[
\lim_{n\to\infty}\|M^{-n}\|^{1/n}=\rho,
\qquad
\lim_{n\to\infty}\|M^n\|^{1/n}=\rho^{-1}.
\]
Since $s_{\min}(M^{-n})=\|M^n\|^{-1}$, we also have
$s_{\min}(M^{-n})^{1/n}\to\rho$.

Choose $k\ge1$ such that $\|M^{-k}\|<1$, and consider
the system indexed by all words $u\in D_N^k$.
Each map $S_u$ has linear part $M^{-k}$ and is therefore
contractive in the Euclidean norm. Iterating the
self-affine identity gives
$
\mu=N^{-k}\sum_{u\in D_N^k}(S_u)_\#\mu.
$
Thus the block system has the same invariant measure $\mu$.
We treat the $N^k$ words as distinct symbols,
even if some corresponding maps coincide. Every composition of $n$ block maps has derivative $M^{-kn}$.
The preceding limits give
\[
-\lim_{n\to\infty}\frac1n\log\|M^{-kn}\|
=
-\lim_{n\to\infty}\frac1n\log s_{\min}(M^{-kn})
=-k\log\rho.
\]
Hence the upper and lower Lyapunov exponents of the block
system coincide and equal $-k\log\rho$ at every coding sequence.

Under the block coding, $\PP$ is a uniform Bernoulli measure
and is ergodic for $\vartheta^k$.
The block entropy calculation above gives
$h_\pi(\vartheta^k,\PP)=k\log N-H_k$.
Applying \cite[Theorem~2.8]{FH} and using
$H_k=kH$ from \eqref{eq:entropy-chain}, we conclude that
$\mu$ is exact dimensional and
\[
\dim_{\rm H}\mu
=\frac{h_\pi(\vartheta^k,\PP)}{-k\log\rho}
=\frac{k\log N-H_k}{-k\log\rho}
=\frac{\log N-H}{-\log\rho}.
\]
Since $-\log\rho>0$, the dimension equals
$\log N/(-\log\rho)$ if and only if $H=0$.
The equivalences now follow from
Proposition~\ref{prop:coding-rigidity}.	
\end{proof}

For $M=\rho^{-1}I_d$,
Proposition~\ref{prop:coding-rigidity} and
Corollary~\ref{cor:dimension-equivalence} show that the three
conclusions of \cite[Theorem~1.1]{Wang26} are equivalent without
assuming the existence of a frame measure. In our notation, these are the
almost-everywhere Dirac property of $\PP_x$ (equivalently, the
Borel bijection between full-measure subsets of $\Omega$ and $K$),
$H_n=0$ for every $n\ge1$, and exact dimensionality with
$\dim_{\rm H}\mu=\log N/(-\log\rho)$, respectively.
For the one-dimensional dimension formula, see
\cite[Remark~3.9]{WT} and \cite[Proposition~4.7]{Wang26}.

For a general expanding matrix, the dimension formula
\eqref{eq:similarity-dimension} need not hold with
$\rho=|\det M|^{-1/d}$; see
Example~\ref{ex:anisotropic-dimension}.
The next proposition gives an entropy bound in terms of
$|\det M|$ that holds for every expanding matrix $M$.
For a countable measurable partition $\mathcal Q$ of $\R^d$,
write
$H_\mu(\mathcal Q)=-\sum_{Q\in\mathcal Q}\mu(Q)\log\mu(Q)$,
with $0\log0=0$.

\begin{proposition}\label{prop:determinant}
Let $\mu=\mu_{M,\mathscr D}$ be as in \eqref{eq:general-ifs}.
Then
\[
H\ge\log N-\log|\det M|.
\]
\end{proposition}

\begin{proof}
First assume that $\|M^{-1}\|<1$, where $\|\cdot\|$ denotes
the Euclidean operator norm. Put $C_j=[0,1)^d+j$ for
$j\in\Z^d$, and let
$\mathcal Q_n=\{M^{-n}C_j:j\in\Z^d\}$.
Since $\mu=\pi_\#\PP$, \cite[Theorem~2.3(i)]{FH} and
\eqref{eq:projection-entropy} give
\[
\log N-H
=h_\pi(\vartheta,\PP)
=\lim_{n\to\infty}\frac1nH_\mu(\mathcal Q_n).
\]
Let $J_n=\{j\in\Z^d:C_j\cap M^nK\ne\varnothing\}$
and $L_n=\#J_n$.
Thus $L_n$ is the number of cells of $\mathcal Q_n$ meeting
$K$, and it is finite because $M^nK$ is bounded.
For $j\in J_n$, choose $z\in C_j\cap M^nK$.
Every $y\in C_j$ satisfies $y-z\in[-1,1]^d$, so
$C_j\subset M^nK+[-1,1]^d$.
The cubes $C_j$, $j\in J_n$, are pairwise disjoint and each
has Lebesgue measure one. Hence
\[
L_n
\le |M^nK+[-1,1]^d|
=|\det M|^n|K+M^{-n}[-1,1]^d|.
\]
Since $K$ is compact and
$\|M^{-n}\|\le\|M^{-1}\|^n\le1$, the sets
$K+M^{-n}[-1,1]^d$ are contained in a fixed ball.
It follows that $L_n\le C_0|\det M|^n$ for some
$C_0>0$ independent of $n$. Since $\mu(K)=1$, at most $L_n$ cells of $\mathcal Q_n$
have positive $\mu$-measure. Therefore
\[
H_\mu(\mathcal Q_n)
\le\log L_n
\le n\log|\det M|+\log C_0.
\]
Dividing by $n$ and letting $n\to\infty$ gives
$\log N-H\le\log|\det M|$, as required.

Now let $M$ be a general expanding matrix.
Since $M^{-n}\to0$, choose $k\ge1$ such that
$\|M^{-k}\|<1$.
The $k$-step system, indexed by all words $u\in D_N^k$,
has $N^k$ labels of equal weight $N^{-k}$, common linear
part $M^{-k}$, and invariant measure $\mu$.
All word labels are kept, even if some maps coincide.
Its first-level conditional entropy is $H_k=kH$ by
Proposition~\ref{prop:coding-rigidity}.
Applying the preceding argument gives
\[
kH=H_k
\ge\log(N^k)-\log|\det(M^k)|
=k\bigl(\log N-\log|\det M|\bigr).
\]
Dividing by $k$ completes the proof.	
\end{proof}

\subsection{Exact overlaps and translations}

The next lemma follows from
Proposition~\ref{prop:coding-rigidity}; see also
\cite[Proposition~3.10]{WT}.

\begin{lemma}\label{lem:exact-overlap}
If the iterated function system has an exact overlap, then $H>0$.	
\end{lemma}

\begin{proof}
Suppose that $S_u=S_v$ for distinct words $u,v\in D_N^n$.
Then $\mu_u=\mu_v=N^{-n}(S_u)_\#\mu$, and both measures
have total mass $N^{-n}>0$.
If $H=0$, Proposition~\ref{prop:coding-rigidity} would imply
$\mu_u\perp\mu_v$, which is impossible because the two measures
are equal and nonzero.
Thus $H>0$.	
\end{proof}

We next introduce the notation used in the Fourier energy
estimate. For a Bessel measure $\Gamma$ for $\mu$ with bound $C$,
write
\[
\mathcal E_\Gamma(f)
=\int_{\R^d}|\widehat{f\mu}(\lambda)|^2\,d\Gamma(\lambda),
\qquad f\in L^2(\mu).
\]
Thus $\mathcal E_\Gamma(f)\le C\|f\|_{L^2(\mu)}^2$. For $t\in\R^d$, let $T_t(x)=x+t$.
If $\theta$ is a finite Borel measure on $\R^d$, the definitions
of pushforward and Fourier transform give
\begin{equation}\label{eq:translation}
\widehat{(T_t)_\#\theta}(\lambda)
=\int_{\R^d}\ee^{-2\pi i\langle\lambda,x+t\rangle}\,d\theta(x)
=\ee^{-2\pi i\langle\lambda,t\rangle}\widehat\theta(\lambda).
\end{equation}
Thus a measure and its translates have Fourier transforms
with the same modulus. We will use this identity to compare their Fourier energies in the argument below.

\section{The entropy obstruction}\label{sec:overlap}

\subsection{A uniformly separating direction}

On the line, for each $n$, a sufficiently small interval near
an endpoint meets only one level-$n$ branch.
To obtain a similar restriction in $\R^d$, we first choose
a fixed vector $v$ for which the numbers
$\langle v,M^jd_i\rangle$, $i\in D_N$, are separated by
a positive amount independent of $j$.

\begin{lemma}\label{lem:uniform-gap}
There are $v\in\R^d\setminus\{0\}$ and $\delta>0$ such that
\begin{equation}\label{eq:uniform-gap}
|\langle v,M^j(d_i-d_k)\rangle|\ge\delta
\qquad(j\ge0,\ i\ne k).
\end{equation}	
\end{lemma}

\begin{proof}
Since $M$ is expanding, the spectral-radius formula gives
$B>0$ and $0<r<1$ such that $\|M^{-j}\|\le Br^j$ for all
$j\ge0$.
For every $z\ne0$, we have
$\|z\|=\|M^{-j}M^jz\|\le Br^j\|M^jz\|$.
Since the digits are distinct, it follows that
\[
L:=\sum_{j=0}^{\infty}\sum_{0\le i<k<N}
\frac1{\|M^j(d_i-d_k)\|}
\le \frac{B}{1-r}\sum_{i<k}\frac1{\|d_i-d_k\|}
<\infty.
\]

Let $V$ be a standard Gaussian vector in $\R^d$,
and let $\Pr$ denote the corresponding Gaussian probability.
For every $z\ne0$, the random variable $\langle V,z\rangle$
is Gaussian with mean zero and variance $\|z\|^2$.
Its density is bounded above by
$(\sqrt{2\pi}\|z\|)^{-1}$, so
\[
\Pr\bigl(|\langle V,z\rangle|<\delta\bigr)
\le \frac{2\delta}{\sqrt{2\pi}\|z\|}
=\sqrt{\frac2\pi}\,\frac{\delta}{\|z\|}
\qquad(\delta>0).
\]

Choose $\delta>0$ such that $\sqrt{2/\pi}\,\delta L<1$.
By the union bound, the probability that
$|\langle V,M^j(d_i-d_k)\rangle|<\delta$
for at least one $j\ge0$ and $i<k$ is at most
$\sqrt{2/\pi}\,\delta L$.
Thus, with positive probability, all the inequalities in
\eqref{eq:uniform-gap} hold simultaneously.
Choose one such vector $v$.
Since $\delta>0$, this vector is nonzero.	
\end{proof}

\begin{lemma}\label{lem:cap-restriction}
There is a Borel set $E\subset K$ with $c:=\mu(E)>0$
such that, for every $n\ge1$, there exist a word
$w^{(n)}\in D_N^n$ and a Borel set $F_n\subset K$ satisfying
\begin{equation}\label{eq:cap-restriction}
\mu|_{F_n}=N^{-n}(S_{w^{(n)}})_\#(\mu|_E),
\qquad \mu(F_n)=cN^{-n}.
\end{equation}	
\end{lemma}

\begin{proof}
Choose $v$ and $\delta$ as in Lemma~\ref{lem:uniform-gap},
and put
\[
a=\min_{y\in K}\langle v,y\rangle,
\qquad E=\{y\in K:\langle v,y\rangle<a+\delta\}.
\]
The set $E$ is Borel. Since $K$ is compact, there is
$y_0\in K$ with $\langle v,y_0\rangle=a$.
Moreover,
$K\cap B(y_0,\delta/(2\|v\|))\subset E$.
Since $K=\supp\mu$, this gives $c:=\mu(E)>0$.

For each $j\ge0$, put
$\beta_j=\min_{i\in D_N}\langle v,M^jd_i\rangle$.
By \eqref{eq:uniform-gap}, this minimum is attained at a
unique index $i_j$, and
$\langle v,M^jd_i\rangle\ge\beta_j+\delta$ whenever $i\ne i_j$.
Fix $n\ge1$, and set
$s_n=\sum_{j=0}^{n-1}\beta_j$ and
$w^{(n)}=(i_{n-1},\ldots,i_0)$.
By \eqref{eq:word-maps}, we have
$
\langle v,b_w\rangle
=\sum_{\ell=1}^n
\langle v,M^{n-\ell}d_{w_\ell}\rangle
$  for $w\in D_N^n$.
For $w=w^{(n)}$, each summand equals its minimum, so
$\langle v,b_{w^{(n)}}\rangle=s_n$.
For any other word, at least one summand is at least
$\delta$ larger than its minimum, while all remaining
summands are at least their minima.
Hence $\langle v,b_w\rangle\ge s_n+\delta$
whenever $w\ne w^{(n)}$. Define the Borel set
\[
F_n=\{x\in K:\langle v,M^nx\rangle<a+s_n+\delta\}.
\]
For $w\in D_N^n$ and $y\in K$, we have $S_w(y)\in K$ and
$
\langle v,M^nS_w(y)\rangle
=\langle v,y\rangle+\langle v,b_w\rangle.
$
If $w=w^{(n)}$, this quantity is less than
$a+s_n+\delta$ exactly when $y\in E$.
If $w\ne w^{(n)}$, it is at least $a+s_n+\delta$,
since $\langle v,y\rangle\ge a$.
Thus $S_w^{-1}(F_n)\cap K$ equals $E$ for $w=w^{(n)}$
and is empty for all other words of length $n$.

Using $\mu=N^{-n}\sum_{|w|=n}(S_w)_\#\mu$, we obtain,
for every Borel set $B\subset\R^d$,
\begin{align*}
\mu(B\cap F_n)
=N^{-n}\sum_{|w|=n}
\mu\bigl(S_w^{-1}(B\cap F_n)\bigr)
=N^{-n}\mu\bigl(E\cap S_{w^{(n)}}^{-1}(B)\bigr).
\end{align*}
This proves the first identity in \eqref{eq:cap-restriction}.
Taking $B=\R^d$ gives $\mu(F_n)=cN^{-n}$.	
\end{proof}

The word $w^{(n)}$ need not be a prefix of $w^{(n+1)}$.
The argument below uses \eqref{eq:cap-restriction} with
the same set $E$ for every $n$.

\subsection{The Fourier energy estimate}

We use the sets from Lemma~\ref{lem:cap-restriction} to
construct unit-norm test functions. Their Fourier energies
are estimated by the translation argument in
\cite[Lemma~3.3]{WT}.

\begin{proposition}
Let $E$, $c$, and $F_n$ be as in
Lemma~\ref{lem:cap-restriction}, and put
$f_n=\ind_{F_n}/\sqrt{\mu(F_n)}$ for $n\ge1$.
If $\Gamma$ is a Bessel measure for $\mu$ with bound $C$,
then, for every $n\ge1$,
\begin{equation}\label{eq:energy-decay}
\|f_n\|_{L^2(\mu)}=1,
\qquad
\mathcal E_\Gamma(f_n)\le\frac Cc R_n.
\end{equation}	
\end{proposition}

\begin{proof}
Fix $n\ge1$. Since $\mu(F_n)=cN^{-n}>0$, we have
$\|f_n\|_{L^2(\mu)}^2
=\mu(F_n)^{-1}\int\ind_{F_n}\,d\mu=1$.
Set $\theta_n=\mu|_{F_n}$.
For each word $w\in D_N^n$, define
$
\theta_{n,w}=N^{-n}(S_w)_\#(\mu|_E).
$
Since $\mu|_E\le\mu$, we have
$0\le\theta_{n,w}\le\mu_w\le\mu$.
Thus the density $h_{n,w}:=d\theta_{n,w}/d\mu$ exists
and satisfies $0\le h_{n,w}\le p_w\le1$ $\mu$-almost
everywhere. In particular, $h_{n,w}\in L^2(\mu)$.
Since $S_w(y)=M^{-n}y+t_w$, we have
$S_w=T_{t_w-t_{w^{(n)}}}\circ S_{w^{(n)}}$.
Together with \eqref{eq:cap-restriction}, this gives
$
\theta_{n,w}
=(T_{t_w-t_{w^{(n)}}})_\#\theta_n.
$
By \eqref{eq:translation}, the Fourier transforms of
$\theta_{n,w}$ and $\theta_n$ have the same modulus.
Since $h_{n,w}\mu=\theta_{n,w}$ and
$f_n\mu=\theta_n/\sqrt{\mu(F_n)}$, it follows that
\[
\mathcal E_\Gamma(h_{n,w})
=\int_{\R^d}|\widehat{\theta_n}(\lambda)|^2
\,d\Gamma(\lambda)
=\mu(F_n)\mathcal E_\Gamma(f_n).
\]

There are $N^n$ words of length $n$.
Summing over these words and applying the Bessel inequality,
we obtain
\begin{align*}
N^n\mu(F_n)\mathcal E_\Gamma(f_n)
=\sum_{|w|=n}\mathcal E_\Gamma(h_{n,w})
\le C\sum_{|w|=n}\int h_{n,w}^2\,d\mu
\le C\int\sum_{|w|=n}p_w^2\,d\mu
=CR_n.
\end{align*}
Since $N^n\mu(F_n)=c$, this proves
\eqref{eq:energy-decay}.	
\end{proof}

\begin{proof}[Proof of Theorem~\ref{thm:general-digits}]
Assume that $H>0$, and suppose that $\Gamma$ is a frame
measure for $\mu$ with bounds $0<A\le C<\infty$.
By Lemma~\ref{lem:prefix-probabilities}, we have $R_n\to0$.
Since $\|f_n\|_{L^2(\mu)}=1$, the lower frame bound and
\eqref{eq:energy-decay} give
\[
0<A\le\mathcal E_\Gamma(f_n)
\le\frac Cc R_n\longrightarrow0,
\]
a contradiction.	
\end{proof}

\begin{proof}[Proof of Corollary~\ref{thm:overlap}]
If $N>|\det M|$, Proposition~\ref{prop:determinant} gives
$H\ge\log N-\log|\det M|>0$.
The conclusion follows from
Theorem~\ref{thm:general-digits}.	
\end{proof}

\section{Critical measures and examples}\label{sec:critical}

\subsection{The critical lattice case}

The equivalence of \textup{(iii)} and \textup{(iv)} in
Theorem~\ref{thm:critical} follows from
\cite[Theorem~1.1 and Corollary~1.1a]{LW}.
The assumption $0\in\mathscr D$ in that reference can be
met by replacing $\mathscr D$ with $\mathscr D-d_0$.
This translates $K$ by $-(M-I_d)^{-1}d_0$ and translates
each finite expansion set, without changing its cardinality
or $|K|$.
We give a direct counting proof of
\textup{(iii)}$\Rightarrow$\textup{(iv)} below.

\begin{proof}[Proof of Theorem~\ref{thm:critical}]
For $n\ge1$, put
$
\mathscr D_n
=\mathscr D+M\mathscr D+\cdots+M^{n-1}\mathscr D.
$
By \eqref{eq:word-maps}, this is the set of vectors $b_w$
with $w\in D_N^n$.
For words $u,v$ of length $n$, the equality $S_u=S_v$
holds if and only if $b_u=b_v$.
Since there are $N^n$ words of length $n$, the absence
of exact overlaps is equivalent to
$\#\mathscr D_n=N^n$ for every $n\ge1$.
This proves \textup{(ii)}$\Longleftrightarrow$\textup{(iii)}.

Assume \textup{(iii)}, and fix $y_0\in K$.
For each $b\in\mathscr D_n$, we have
$M^{-n}(y_0+b)\in K$, so $b\in M^nK-y_0$.
Moreover, $\mathscr D_n\subset\Z^d$ has $N^n$ elements.
The half-open cubes $b+[0,1)^d$, $b\in\mathscr D_n$,
are therefore pairwise disjoint, each has volume one,
and all are contained in $M^nK-y_0+[0,1]^d$.
Hence
\[
N^n
\le |M^nK-y_0+[0,1]^d|
=|\det M|^n\,|K+M^{-n}[0,1]^d|,
\]
where the equality follows from a linear change of variables
and translation invariance of Lebesgue measure.
Since $N=|\det M|$, we obtain
$|K+M^{-n}[0,1]^d|\ge1$. For every $\varepsilon>0$, the relation $\|M^{-n}\|\to0$
gives
$M^{-n}[0,1]^d\subset[-\varepsilon,\varepsilon]^d$
for all sufficiently large $n$.
Thus $|K+[-\varepsilon,\varepsilon]^d|\ge1$.
Since $K$ is compact, these sets decrease to $K$ as
$\varepsilon\downarrow0$.
Continuity from above gives $|K|\ge1$, proving
\textup{(iv)}.

Assume \textup{(iv)}.
Since $K=\bigcup_{i=0}^{N-1}S_i(K)$ and
$
\sum_{i=0}^{N-1}|S_i(K)|
=N|\det M|^{-1}|K|
=|K|,
$
the sets $S_i(K)$ overlap only on sets of Lebesgue measure
zero.
Put $\nu=|K|^{-1}\ind_K\,dx$.
For every Borel set $B\subset\R^d$, a change of variables
gives
\begin{align*}
\frac1N\sum_{i=0}^{N-1}(S_i)_\#\nu(B)
=\frac1{N|K|}
\sum_{i=0}^{N-1}|S_i^{-1}(B)\cap K|
=\frac{|\det M|}{N|K|}
\sum_{i=0}^{N-1}|B\cap S_i(K)|
=\frac{|B\cap K|}{|K|}
=\nu(B).
\end{align*}
Thus $\nu$ satisfies \eqref{eq:general-ifs}.
By uniqueness of the invariant probability measure,
$\mu=\nu$. This proves \textup{(v)}.

Assume \textup{(v)}, and let $f\in L^2(\mu)$.
Extend $f$ by zero outside $K$.
Then $f\mu=|K|^{-1}f\ind_K\,dx$, and Plancherel's theorem
gives
\[
\int_{\R^d}|\widehat{f\mu}(\lambda)|^2
|K|\,d\lambda
=\frac1{|K|}\int_K|f(x)|^2\,dx
=\|f\|_{L^2(\mu)}^2.
\]
Hence \textup{(v)} implies \textup{(i)}, with the Parseval
frame measure $d\Gamma(\lambda)=|K|\,d\lambda$.

Finally, Theorem~\ref{thm:general-digits} gives
\textup{(i)}$\Rightarrow$\textup{(vi)}, and
Lemma~\ref{lem:exact-overlap} gives
\textup{(vi)}$\Rightarrow$\textup{(ii)}.
All six conditions are therefore equivalent.	
\end{proof}

\subsection{Dimension and the matrix hypothesis}

Theorem~\ref{thm:general-digits} and
Corollary~\ref{cor:dimension-equivalence} give the following result.

\begin{corollary}\label{cor:similarities}
Suppose that every eigenvalue of $M$ has modulus
$\rho^{-1}>1$.
If $\mu_{M,\mathscr D}$ admits a frame measure, then
\[
\dim_{\rm H}\mu_{M,\mathscr D}
=\frac{\log N}{-\log\rho}\le d.
\]
Consequently, $\mu_{M,\mathscr D}$ admits no frame measure
if $N\rho^d>1$.

If $M=\rho^{-1}I_d$, put
\[
V=\operatorname{span}\{d_i-d_0:1\le i<N\},
\qquad k=\dim V.
\]
The affine hull $\operatorname{aff}\mathscr D=d_0+V$
is the smallest affine subspace of $\R^d$ containing
$\mathscr D$, and its dimension is $k$.
In this case, $\mu_{M,\mathscr D}$ admits no frame measure
if $N\rho^k>1$.	
\end{corollary}

\begin{proof}
Write $\mu=\mu_{M,\mathscr D}$, and suppose that $\mu$
admits a frame measure.
By Theorem~\ref{thm:general-digits}, we have $H=0$.
Thus \eqref{eq:similarity-dimension} gives
\[
\dim_{\rm H}\mu
=\frac{\log N-H}{-\log\rho}
=\frac{\log N}{-\log\rho}.
\]
Since $\mu$ is supported on $\R^d$, we have
$\dim_{\rm H}\mu\le d$.
Hence $\log N/(-\log\rho)\le d$, or equivalently
$N\rho^d\le1$.
This proves the first two assertions.

Now suppose that $M=\rho^{-1}I_d$, and let $V$ be as above.
Put $x_0=\rho d_0/(1-\rho)$.
For every $\omega\in\Omega$, the coding formula gives
\[
\pi(\omega)
=\sum_{j\ge1}\rho^j d_{\omega_j}
=x_0+\sum_{j\ge1}\rho^j(d_{\omega_j}-d_0).
\]
Every term of the last series lies in $V$.
Since $V$ is closed, its sum also belongs to $V$.
Therefore $K=\pi(\Omega)\subset x_0+V$, and
$\dim_{\rm H}\mu\le\dim_{\rm H}K\le k$.
By \eqref{eq:similarity-dimension}, we obtain
$\log N/(-\log\rho)\le k$, or equivalently
$N\rho^k\le1$.
Thus $N\rho^k>1$ rules out frame measures.	
\end{proof}

The dimension conclusion and the nonexistence condition
$N\rho^d>1$ apply, in particular, when $M^{-1}=\rho O$
with $O$ orthogonal.
For $O=I_d$, these conclusions were proved in
\cite[Theorem~1.1 and Corollary~1.5]{Wang26}.
For a general $O$, the space spanned by the digit differences
need not be invariant under $O$, so $K$ need not lie
in a $k$-dimensional affine subspace.

If $O^r=I_d$ for some $r\ge1$ and there is no exact overlap,
the $r$-step system consists of $N^r$ distinct maps
$S_w(x)=\rho^r x+t_w$, $|w|=r$, each with weight $N^{-r}$.
It has the same invariant measure $\mu$, and its first-block
conditional entropy is $H_r=rH$ by \eqref{eq:entropy-chain}.
Applying \cite[Theorem~1.1]{Wang26} to this system shows
that frame-measure existence forces $H_r=0$, and hence $H=0$.
The proof of Theorem~\ref{thm:general-digits} does not
require any power of $M$ to be scalar.

\begin{example}\label{ex:anisotropic-dimension}
Let
\[
M=\begin{pmatrix}2&0\\0&4\end{pmatrix},
\qquad
\mathscr D=\{0,1\}\times\{0,2\}.
\]
Here $N=4$ and $|\det M|=8$.
The product structure of the digit set gives
$
\mu=(\ind_{[0,1]}\,dx)\otimes\nu,
$
where $\nu$ is the uniform Cantor measure generated by
$y\mapsto y/4$ and $y\mapsto(y+2)/4$.
Put $C=\supp\nu$. Then $C\subset[0,2/3]$, and its two
first-level pieces lie in the disjoint intervals
$[0,1/6]$ and $[1/2,2/3]$.
The two horizontal intervals are $[0,1/2]$ and $[1/2,1]$;
their intersection has Lebesgue measure zero.
Thus the four first-level component measures of $\mu$
are pairwise mutually singular, and
Proposition~\ref{prop:coding-rigidity} gives $H=0$.

Applying \eqref{eq:similarity-dimension} to the two-map
system defining $\nu$, whose conditional entropy is zero,
shows that $\nu$ is exact dimensional with
$\dim_{\rm H}\nu=\log2/\log4=1/2$.
The first factor is exact dimensional with dimension one.
To compute the dimension of their product, write
$I_r(t)=(t-r,t+r)$ and note that
\[
I_{r/\sqrt2}(x)\times I_{r/\sqrt2}(y)
\subset B((x,y),r)
\subset I_r(x)\times I_r(y).
\]
For $\mu$-almost every $(x,y)$, the logarithm of the measure
of either rectangle, divided by $\log r$, tends to $1+1/2$
as $r\downarrow0$. These inclusions therefore show that $\mu$ is
exact dimensional and
$
\dim_{\rm H}\mu=1+\frac 12=\frac32.
$

On the other hand, setting
$\rho_*=|\det M|^{-1/2}=8^{-1/2}$ gives
$
\frac{\log N}{-\log\rho_*}=\frac43\ne\frac32.
$
Moreover, $K=[0,1]\times C$ contains the three noncollinear
points $(0,0)$, $(1,0)$, and $(0,2/3)$, so
$\operatorname{aff}K=\R^2$.
Thus, even when $H=0$ and the attractor has full affine span,
one cannot extend \eqref{eq:similarity-dimension} to arbitrary
expanding matrices by setting $\rho=|\det M|^{-1/d}$.	
\end{example}

\begin{example}\label{ex:zero-not-frame}
Let $M=3$ and $\mathscr D=\{0,1\}$.
The attractor is contained in $[0,1/2]$, and
\[
S_0(K)\subset[0,1/6],
\qquad
S_1(K)\subset[1/3,1/2].
\]
The first-level pieces are therefore disjoint.
Hence $H=0$, and all conditions in
Proposition~\ref{prop:coding-rigidity} hold.
Corollary~\ref{cor:dimension-equivalence} also gives
$\dim_{\rm H}\mu=\log2/\log3$. However, the integer-power obstruction of Fu and Wang
\cite[Theorem~1.1]{FW} applies with $N=2$, $\rho=1/3$, and $m=1$,
since $\rho^{-1}=3$ and $2\nmid3$.
Thus $\mu$ admits no frame measure.
This shows that the zero-entropy conditions, including
the corresponding dimension equality, are not sufficient
for frame-measure existence.	
\end{example}

\begin{example}\label{ex:jordan}
Let
\[
q=\frac\pi2,\qquad
B=\begin{pmatrix}1&1\\0&1\end{pmatrix},\qquad
M=qB,\qquad
\mathscr D=\{(0,0),(1,0),(0,1)\}.
\]
Both eigenvalues of $M$ equal $q>1$, and $M$ is not
diagonalizable.
Since $N=3>|\det M|=q^2$, Proposition~\ref{prop:determinant}
gives $H\ge\log(3/q^2)>0$, and
Corollary~\ref{thm:overlap} rules out frame measures. There is no exact overlap.
Indeed, suppose that $S_u=S_v$ for words $u,v$ of length $n$.
Then $b_u=b_v$, so
\[
\sum_{j=0}^{n-1}
q^jB^j(d_{u_{n-j}}-d_{v_{n-j}})=0.
\]
Each coordinate is an integer polynomial in $q$.
Since $q$ is transcendental, every coefficient is zero.
Thus
$B^j(d_{u_{n-j}}-d_{v_{n-j}})=0$ for every $0\le j<n$.
The invertibility of $B^j$ gives $d_{u_{n-j}}=d_{v_{n-j}}$
for every $j$, and hence $u=v$.
Thus positive conditional entropy can occur without exact
overlaps, also for non-diagonalizable matrices.
Moreover, $M^r=q^r\bigl(\begin{smallmatrix}1&r\\0&1\end{smallmatrix}\bigr)$
is non-diagonalizable for every $r\ge1$, so no power of $M$
is similar to a scalar matrix.
\end{example}

\begin{example}\label{ex:nonstandard-digits}
Let $M=2I_2$ and
$\mathscr D=\{0,2\}\times\{0,1\}$.
Then $N=4=|\det M|$.
The four images of $[0,2]\times[0,1]$ cover this rectangle
and have disjoint interiors. Consequently,
\[
K=[0,2]\times[0,1],
\qquad
\mu=\tfrac12\ind_K\,dx.
\]

Modulo $2\Z^2$, the four digits occupy only the two classes
represented by $(0,0)$ and $(0,1)$.
Thus $\mathscr D$ is not a complete set of representatives
for $\Z^2/2\Z^2$.
Nevertheless, the Fourier bases on the two intervals give
the spectrum	
$
\Lambda=(\tfrac12\Z)\times\Z,	
$
and Plancherel's theorem gives the Parseval frame measure
$d\Gamma(\lambda)=2\,d\lambda$.
Hence the complete residue condition is unnecessary in
Theorem~\ref{thm:critical}.

This example also has $\rho=1/2>1/N=1/4$, while
$N\rho^2=1$.
Since $\mu$ is spectral, the one-dimensional condition
$\rho>1/N$ cannot be used unchanged in higher dimensions.	
\end{example}

\subsection{The common-matrix hypothesis}

The Fourier energy estimate uses the fact that, at each
level $n$, the measures $\theta_{n,w}$ are translates of
one another. This follows from the common linear part
$M^{-n}$ of the maps $S_w$.
The following example shows that positive conditional
entropy need not exclude a frame measure when the
contraction ratios vary, even if all weights are equal.

\begin{example}\label{ex:variable-contractions}
Consider the five maps on $\R$
\[
S_0(x)=\frac x2,\quad
S_1(x)=\frac{x+1}2,\quad
S_2(x)=\frac x3,\quad
S_3(x)=\frac{x+1}3,\quad
S_4(x)=\frac{x+2}3,
\]
each with weight $1/5$.
Let $X$ be the coding random variable for this system.

Put $\mu=\ind_{[0,1]}\,dx$.
The first two pushforwards have density $2$ on the two
half intervals, and their densities sum to
$2\ind_{[0,1]}$ almost everywhere.
Similarly, the densities of the last three pushforwards
sum to $3\ind_{[0,1]}$.
Therefore
$
\frac15\sum_{i=0}^4(S_i)_\#\mu
=\frac15(2+3)\mu
=\mu.
$
Thus $\mu$ is the invariant probability measure.

The conditional probabilities are the weighted branch
densities
$p_i=\frac15\,d((S_i)_\#\mu)/d\mu$.
For $\mu$-almost every $x$, exactly one half-interval branch
and one third-interval branch contribute at $x$.
The conditional probability vector therefore has precisely
two nonzero entries, $2/5$ and $3/5$.
Hence
\[
H(\omega_1\mid X)
=-\frac25\log\frac25-\frac35\log\frac35>0.
\]
Nevertheless, $\Z$ is a spectrum for $\mu$, so
$\sum_{m\in\Z}\delta_m$ is a Parseval frame measure.	
\end{example}

For systems with a common expanding matrix, Corollary~\ref{thm:overlap} rules out frame measures when $N>|\det M|$, and Theorem~\ref{thm:critical} characterizes their existence in the critical integer case. When $N<|\det M|$ and $H=0$, further conditions are needed to decide
frame-measure existence, as
Example~\ref{ex:zero-not-frame} illustrates.
The critical case with arbitrary real digits is not
settled by the lattice-counting argument used here.


\begin{thebibliography}{99}

\bibitem{DC}
Q.-R. Deng and J.-B. Chen,
Uniformity of spectral self-affine measures,
\emph{Adv. Math.} \textbf{380} (2021), Paper No.~107568.

\bibitem{DHW}
D. E. Dutkay, D. Han and E. Weber,
Continuous and discrete Fourier frames for fractal measures,
\emph{Trans. Amer. Math. Soc.} \textbf{366} (2014), 1213--1235.

\bibitem{DL}
D. E. Dutkay and C.-K. Lai,
Uniformity of measures with Fourier frames,
\emph{Adv. Math.} \textbf{252} (2014), 684--707.

\bibitem{FH}
D.-J. Feng and H. Hu, Dimension theory of iterated function systems,
\emph{Comm. Pure Appl. Math.} \textbf{62} (2009), 1435--1500.

\bibitem{FW}
X.-Y. Fu and W.-J. Wang,
Nonexistence of frame measures for separated uniform
consecutive-digit Bernoulli convolutions,
preprint, \mbox{\href{https://arxiv.org/abs/2608.19676}
	{arXiv:2608.19676}}, 2026.

\bibitem{Hutchinson}
J. E. Hutchinson, Fractals and self-similarity,
\emph{Indiana Univ. Math. J.} \textbf{30} (1981), 713--747.

\bibitem{Lai11}
C.-K. Lai, On Fourier frame of absolutely continuous measures,
\emph{J. Funct. Anal.} \textbf{261} (2011), 2877--2889.

\bibitem{LW}
J. C. Lagarias and Y. Wang, Self-affine tiles in $\R^n$,
\emph{Adv. Math.} \textbf{121} (1996), 21--49.

\bibitem{Wang26}
W.-J. Wang,
Frame measure rigidity and fiber entropy for homogeneous self-similar measures,
preprint, \mbox{\href{https://arxiv.org/abs/2611.14002}{arXiv:2611.14002}}, 2026.

\bibitem{WT}
Z. Wang and W. Tang, Entropy and frame measures for self-similar measures,
preprint, \mbox{\href{https://arxiv.org/abs/2611.09190}{arXiv:2611.09190}}, 2026.

\end{thebibliography}
\end{document}